\documentclass[a4,twoside,10pt]{amsart}

\usepackage{amssymb}
\usepackage{amsthm}
\usepackage[english]{babel}
\usepackage{enumerate}
\usepackage{fontenc}
\usepackage[colorlinks=true,linkcolor=blue,citecolor=red,urlcolor=blue]{hyperref}
\usepackage[latin1]{inputenc}
\usepackage{mathrsfs}
\usepackage{orcidlink}
\usepackage{physics}
\usepackage{xcolor}
\usepackage[english]{varioref}
\usepackage[all]{xy}
\usepackage{url}

\theoremstyle{plain}
\newtheorem{theorem}{Theorem}[section]
\newtheorem{conjecture}[theorem]{Conjecture}
\newtheorem{corollary}[theorem]{Corollary}
\newtheorem{lemma}[theorem]{Lemma}
\newtheorem{proposition}[theorem]{Proposition}

\theoremstyle{definition}
\newtheorem{definition}[theorem]{Definition}
\newtheorem{esempio}[theorem]{Example}
\newenvironment{example}{\begin{esempio}}{\hfill{$\triangle$}\end{esempio}}
\newtheorem{construction}[theorem]{Construction}
\newtheorem*{Proof}{Proof}
\newenvironment{prf}{\begin{Proof}}{\hfill{$\Box$}\end{Proof}}
\newtheorem{rmk}[theorem]{Remark}
\newenvironment{remark}{\begin{rmk}}{\hfill{$\Diamond$}\end{rmk}}
\newtheorem{remarks}[theorem]{Remarks}

\DeclareMathOperator{\Ann}{Ann}

\DeclareMathOperator{\Hilb}{Hilb}
\let\Im\relax
\DeclareMathOperator{\Im}{Im}
\DeclareMathOperator{\pr}{pr}
\DeclareMathOperator{\uzero}{\underline{0}}
\DeclareMathOperator{\wl}{\widetilde{\ell}}
\DeclareMathOperator{\whvarphi}{\widehat{\varphi}}
\DeclareMathOperator{\wtvarphi}{\widetilde{\varphi}}
\DeclareMathOperator{\whpsi}{\widehat{\psi}}
\DeclareMathOperator{\wtpsi}{\widetilde{\psi}}

\newcommand{\preceqlex}{\preceq_{lex}}

\newcommand{\A}{\mathbb{A}}
\newcommand{\B}{\mathbb{B}}
\newcommand{\D}{\mathbb{D}}
\newcommand{\K}{\mathbb{K}}
\newcommand{\N}{\mathbb{N}}
\newcommand{\bP}{\mathbb{P}}
\newcommand{\R}{\mathbb{R}}
\newcommand{\bS}{\mathbb{S}}
\newcommand{\V}{\mathbb{V}}
\newcommand{\W}{\mathbb{W}}
\newcommand{\Z}{\mathbb{Z}}

\newcommand{\cAG}{\mathcal{AG}}

\newcommand{\sP}{\mathscr{P}}

\newcommand{\wA}{\widetilde{A}}
\newcommand{\wG}{\widetilde{G}}
\newcommand{\wH}{\widetilde{H}}
\newcommand{\wh}{\widetilde{h}}
\newcommand{\wL}{\widetilde{L}}
\newcommand{\wP}{\widetilde{P}}

\newcommand{\wx}{\widetilde{x}}

\newcommand{\ua}{\underline{a}}
\newcommand{\ub}{\underline{b}}
\newcommand{\uc}{\underline{c}}
\newcommand{\ud}{\underline{d}}
\newcommand{\ui}{\underline{i}}
\newcommand{\uj}{\underline{j}}
\newcommand{\ul}{\underline{\ell}}
\newcommand{\uv}{\underline{v}}

\begin{document}

\title{CW-complexes and minimal Hilbert vector of standard graded Artinian Gorenstein algebras}

\author{Armando Capasso}
\address{Universit\`a degli Studi di Trieste, P.le Europa 1, Trieste (Italy), C.A.P. 34127}
\email{armando.capasso@units.it}

\thanks{A.C. is member of INdAM - GNSAGA. ORCID: 0009-0001-5463-7221 \orcidlink{0009-0001-5463-7221}}



\date{\today}

\begin{abstract}
I introduce a geometric interpretation of the set of \emph{standard graded Artinian Gorenstein algebras of codimension $n$ and degree $d$}: the \emph{standard locus}, which is a subset of the projective space of degree $d$ polynomials in $n$ variables, and I characterize it. Under opportune hypothesis, I prove that the \emph{locus of full Perazzo polynomials} is the union of the minimal dimensional irreducible components of the standard locus and it is a pure dimensional subset. On the other hand, I associate to any homogeneous polynomial a topological space, which is a CW-complex. Using all these sets, I prove that the Hilbert function restricted to the standard locus has minimal values on any irreducible component of the domain. I apply all this to the \emph{Full Perazzo Conjecture} and I prove it.
\end{abstract}

\maketitle

\tableofcontents

\section*{Introduction}
\markboth{Introduction}{Introduction}

\noindent Poincar\'e algebras arise as de Rham cohomology algebras of compact oriented manifolds without boundary. These are graded Artinian algebras, and where the underlying manifold is connected, these algebras are also Gorenstein. \emph{Vice versa}, a graded Artinian Gorenstein algebra (\emph{GAG} algebras, for short) is a Poincar\'e algebra (\cite[Proposition 2.79]{HMMNWW}). I deal with standard graded Artinian Gorenstein $\K$-algebras (\emph{SGAG} algebras, for short) over a field of characteristic zero.\smallskip

\noindent A natural and classical problem consists in understanding their possible Hilbert functions, sometimes also called Hilbert vectors. By duality, these vectors are \emph{symmetric}, \emph{i.e.} let $\left(h_0,h_1,\dotsc,h_d\right)$ be a Hilbert vector of a codimension $n$ SGAG algebra of socle degree $d$, then $h_i=h_{d-i}$ for any $i\in\{0,\dotsc,d\}$, $h_0=h_d=1$ and $h_1=n$. When the codimension of the SGAG algebra is less than or equal to $3$, all possible Hilbert vectors were characterized by \cite[Theorem 4.2]{SR:2}, and these are all \emph{unimodal} vectors, \emph{i.e.} a symmetric vector $\left(h_0,h_1,\dotsc,h_d\right)$ is unimodal if $h_0\leq h_1\leq\dotsc\leq h_t$ where $\displaystyle t=\left\lfloor\frac{d}{2}\right\rfloor$. On the other hand, it is known that non-unimodal vectors exist for any SGAG algebra of codimension greater than or equal to $5$ (\cite[Section 2]{B:I} and \cite[Proposition 5.3 and Remark 5.4]{B:L}). It is an open problem whether non unimodal Hilbert vectors of codimension $4$ SGAG algebras exist.\smallskip

\noindent Let $\cAG_{st,\K}(n,d)$ be the set of SGAG $\K$-algebras of socle degree $d$ and codimension $n$. There is a natural partial order between the Hilbert vectors of SGAG algebras: one gets
\begin{displaymath}
\left(1,n,h_2,\dotsc,h_2,n,1\right)\preceq\left(1,n,\wh_2,\dotsc,\wh_2,n,1\right)
\end{displaymath}
if $h_i\leq\wh_i$ for any $\displaystyle i\in\left\{0,\dotsc,\left\lfloor\frac{d}{2}\right\rfloor\right\}$. This order makes $\cAG_{st,\K}(n,d)$ a \emph{poset}, \emph{i.e.} $A\preceq\wA$ if and only if $\Hilb(A)\preceq\Hilb\left(\wA\right)$, where the notations are obvious. Thus it is natural to look for the minimal and maximal elements of $\cAG_{st,\K}(n,d)$. On one hand, maximal SGAG algebras are characterized by \cite[Definition 3.11 and Proposition 3.12]{I:K}. On the other hand, in \cite{BGIZ} the authors have conjectured that particular SGAG algebras, called \emph{full Perazzo algebras}, have minimal Hilbert vector (\cite[Conjecture 2.6]{BGIZ}).\smallskip

\noindent To be clear, a \emph{full Perazzo polynomial} is a bidegree $(1,d-1)$ polynomial $\displaystyle {f=\sum_{i=1}^{\tau(n,d-1)}x_iM_i}$ where $\tau(n,d)=\displaystyle\binom{n+d-1}{d}$ and $\left\{M_1,\dotsc,M_{\tau(n,d-1)}\in \K\left[u_1,\dotsc,u_n\right]_{d-1}\right\}$ is the monic monomial basis. The associated SGAG $\K$-algebra is a full Perazzo algebra (see Theorem \ref{th1.1}). It is known that socle degrees $4$, $5$ partially and where $n=3$ full Perazzo algebras are minimal in $\left(\cAG_{st,\K}\left(n+\tau(n,d-1),d\right),\preceq\right)$ (see \cite[Theorem 4.2]{CGIZ} and \cite[Theorems 3.15 and 4.2]{BGIZ}).\smallskip

\noindent In this paper I give a complete proof of \emph{Full Perazzo Conjecture} (Theorem \ref{th5.1}), which is the main result. The relevant proof is based on three bedrocks: the introduction of the so-called \emph{standard locus} (see Equation \eqref{eq1}) and \emph{full Perazzo locus} (see Equation \eqref{eq5}), the association of a \emph{CW-complex} to a homogeneous polynomial (Construction \ref{Constr1}).\smallskip

\noindent The idea behind the \emph{standard locus} is the following one: let $f\in\K\left[x_1,\dotsc,x_m\right]_d$, it is known for any $\lambda\in\K\setminus\{0\}$, $\Ann(f)=\Ann(\lambda f)$ (see Definition \ref{def1.2}) hence the Hilbert vector of the GAG $\K$-algebra determined by $f$ depends only on ${[f]\in\bP\left(\K\left[x_1,\dotsc,x_m\right]_d\right)}$. Thus it is interested to study the set $U^{(m,d)}_{st}$ of all degree $d$ polynomials in $m$ variables such that the relevant GAG $\K$-algebra is also standard, up to multiplication by a non-zero scalar. One carries on this set the Zariski topology and proves its locally closeness (Lemma \ref{lem2.1}), smoothness (Corollary \ref{cor2.1}) and determines its irreducible components (Lemma \ref{lem2.3}). Since $U^{(m,d)}_{st}$ has finitely many irreducible components, in this paper I am interested to whose have minimal dimension. Under the opportune hypothesis, the \emph{full Perazzo locus} $U^{(n,d)}_{FP}$, \emph{i.e.} the set of full Perazzo polynomials up to multiplication by a non-zero scalar, is in the union of the minimal dimensional irreducible components of $U^{(n+\tau(n,d-1),d)}_{st}$ and it is a pure dimensional subset (Lemma \ref{lem2.2}).\smallskip

\noindent Before to explain the ideas behind Construction \ref{Constr1}, I recall that Reisner and Stanley associated a \emph{simplicial complex} to square-free homogeneous polynomials (\cite[Section 3]{SR:1} and \cite[Section 1]{R:GA}), and this construction was partially generalized by Cerminara, Gondim, Ilardi and Maddaloni to bihomogeneous polynomials (\cite[Section 3.1]{CGIM}).\smallskip

\noindent More precisely, I introduce a new Construction \ref{Constr1} which allows me to identify each (monic) monomial of degree $d$ in $m$ variables with an element of the $(d-1)$-skeleton of a CW-complex that I call $P(m)$. To start, consider a monomial $x_1^{k_1}\cdot\dotsc\cdot x_m^{k_m}$, where $k_1+\dotsc+k_m=d$, as a product of $d$ linear forms $\wx_{1,1}\cdot\dotsc\cdot\wx_{m,k_m}$; all $\wx_{j,k_j}$, with fixed index $j$, are identified to a single point $x_j$, one obtains a bouquet of $k_j$-circles with marked point $x_j$ and $\displaystyle\sum_{j^{\prime}\in\{1,\dotsc,m\}\setminus\{j\}}k_jk_{j^{\prime}}\,1$-cells attached to $\wx_j$ and $\wx_{j^{\prime}}$, one identifies all these circles between them and all these $1$-cells between them, and so on (see also \cite[Construction 3.10]{CDPI}). At the end one has a CW-complex $\zeta_{x_1^{k_1}\cdot\dotsc\cdot x_m^{k_m}}$, and attaching all these CW-complexes along the common skeletons one obtains $P(m)$.\smallskip

\noindent Moreover, this construction corrects \cite[Construction 3.10]{CDPI}, since the authors forgot to identify the CW-complexes with skeletons as described above. However, this correction does not invalidate \cite[Theorems 3.16 and 3.18]{CDPI}.\smallskip

\noindent I am not able to give a complete characterization of the annihilator ideal of a homogeneous polynomial and of Hilbert vector of the corresponding SGAG $\K$-algebra, studying the CW-complex associated, which are the main results of \cite{CDPI}. However the key result states that the function
\begin{displaymath}
h\colon[f]\in U^{(n+\tau(n,d-1),d)}_{st}\to\Hilb\left(\K\left[X_1,\dotsc,X_{\tau(n,d-1)},U_1,\dotsc,U_n\right]/\Ann(f)\right)\in\left(\N_{\geq1}\right)^{d+1},
\end{displaymath}
which can be interpreted as the \emph{Hilbert function} restricted to ${\cAG_{st,\K}(n+\tau(n,d-1),d)}$, has minimal values on any irreducible component of $U^{(n+\tau(n,d-1),d)}_{st}$ (lemma \ref{lem5.4}); where one considers the order relation $\preceq$ on the codomain of course.\smallskip

\noindent Thus, using the previous function $h$, I prove that the entries of the Hilbert vector of SGAG $\K$-algebras is minimal on the minimal dimensional irreducible components of $U^{(n+\tau(n,d-1),d)}_{st}$ which correspond to full Perazzo polynomials. All this, under the opportune hypothesis, allows me to achieve the main result of this paper.\smallskip

\noindent\textbf{Notations.} In this the paper I fix the following notations and assumptions, unless otherwise indicated.
\begin{itemize}
\item For any set $Y$, $\mathscr{P}(Y)$ is the set of subsets of $Y$;
\item For any finite set $Y$, its cardinality will be indicated as $\# Y$.
\item For any $k,n\in\N_{\geq0}$, $\displaystyle\binom{n}{k}=\begin{cases}
\displaystyle\frac{n!}{k!(n-k)!}\iff k\leq n\\
0\iff k>n
\end{cases}$.
\item $\K$ is an algebraically closed field of characteristic $0$.
\item $R=\K\left[x_1,\dotsc,x_m\right]$ will always be the graded ring of polynomials in $m$ variables $x_1,\dotsc,x_m$; each $x_i$ has degree $1$.
\item $Q=\K\left[X_1,\dotsc,X_m\right]$ will be the the ring of differential operators of $R$, where $\displaystyle X_i=\pdv{x_i}$.
\item For any finite set $Y$ of either $R$ or $Q$, $\langle Y\rangle$ is the $\K$-vector space generated by $Y$ and $(Y)$ is the ideal generated by $Y$.
\item $S=R\otimes_{\K}\K\left[u_1,\dotsc,u_n\right]= \K\left[x_1,\dotsc,x_m,u_1,\dotsc,u_n\right]$ is the bigraded ring of polynomials in $m+n$ variables $x_1,\dotsc,x_m,u_1,\dotsc,u_n$; each $x_i$ has bidegree $(1,0)$ and $u_j$ has bidegree $(0,1)$.
\item Define $S_{(d_1,d_2)}$ to be the $\K$-vector space of bihomogeneous polynomials $f$ of bidegree $(d_1,d_2)$, \emph{i.e.} $f$ can be written as $\displaystyle\sum_{i=0}^pa_i b_i$, where ${a_i\in\K\left[x_1,\dotsc,x_m\right]_{d_1}}$ and $b_i\in\K\left[u_1,\dotsc,u_n\right]_{d_2}$.
\item $T=Q\otimes_{\K}\K\left[U_1,\dotsc,U_n\right]=\K\left[X_1,\dotsc,X_m,U_1,\dotsc,U_n\right]$ is the (bigraded) ring of differential operators of $S$, where $\displaystyle X_i=\pdv{x_i}$ and $\displaystyle U_j=\pdv{u_j}$; $X_i$ has bidegree $(1,0)$ and $U_j$ has bidegree $(0,1)$.
\item The subscript of a graded $\K$-algebra or homogeneous ideal will indicate the part of that degree; $R_d$ is the $\K$-vector space of the homogeneous polynomials of degree $d$, $Q_{\delta}$ the $\K$-vector space of differential operators of degree $\delta$, and $J_{\ell}$ is the $\ell$-degree part of the ideal $J$.
\end{itemize}
\noindent \textit{Acknowledgments.} I am grateful to Maurizio Brunetti, Pietro De Poi, Rodrigo Gondim, Giovanna Ilardi and Joachim Jelisiejew for their useful remarks.

\section{Graded Artinian Gorenstein algebras}
\markboth{Graded Artinian Gorenstein algebras}{Graded Artinian Gorenstein algebras}

\subsection{Graded Artinian Gorenstein algebras are Poincar\'e algebras}

\begin{definition}
Let $I$ be a homogeneous ideal of $R$ such that $\displaystyle A=R/I=\bigoplus_{i=0}^d A_i$ is a graded Artinian $\K$-algebra, where $A_d\neq 0$, $A_e=0$ for any $e<0$ and $e>d$. The integer $d$ is the \emph{socle degree of $A$}. Setting $h_i=\dim_{\K}A_i$, the vector ${\Hilb(A)=(1,h_1,\dotsc,h_d)}$ is called \emph{Hilbert vector of $A$}. The algebra $A$ is said \emph{standard} if it is generated in degree $1$, \emph{i.e.} if $I_1=0$. If $A$ is standard then $h_1=m$ is called \emph{codimension of $A$}.
\end{definition}
\begin{definition}[{\cite[Definition 2.1]{M:W}}]
A graded Artinian $\K$-algebra $\displaystyle A=\bigoplus_{i=0}^d A_i$ is a \emph{Poincar\'e algebra} if $\dim A_d=1$ and $\cdot\colon A_i\times A_{d-i}\to A_d\cong\K$ is a \emph{perfect pairing} for each $i\in\{0,\dotsc,d\}$.
\end{definition}
\begin{remark}
The Hilbert vector of a Poincar\'e algebra $A$ is \emph{symmetric with respect to $h_{\textstyle\left\lfloor\frac{d}{2}\right\rfloor}$}, \emph{i.e.} $\Hilb(A)=(1,h_1,h_2,\dotsc,h_2,h_1,1)$.
\end{remark}
\noindent I also recall the following definition.
\begin{definition}[{\cite[Theorem 2.79]{HMMNWW}}]\label{def1.1}
A graded Artinian $\K$-algebra $A$ is \emph{Gorenstein} if (and only if) it is a Poincar\'e algebra.
\end{definition}
\noindent For sake of simplicity, a graded Artinian Gorenstein algebra will be called \emph{GAG algebra}.

\subsection{Graded Artinian Gorenstein quotient algebras of $Q$}

For any $d\geq\delta\geq0$ there exists a natural $\K$-bilinear map $B\colon R_d\times Q_{\delta}\to R_{d-\delta}$ defined by differentiation
\begin{displaymath}
B(f,\alpha)=\alpha(f).
\end{displaymath}
\begin{definition}\label{def1.2}
Let $W=\left\langle f_1,\dotsc,f_{\ell}\right\rangle$ be a finite dimensional $\K$-vector subspace of $R$, where $f_1,\dotsc,f_{\ell}$ are forms in $R$. The \emph{annihilator of $W$ in $Q$} is the following homogeneous ideal
\begin{displaymath}
\Ann(W)=\{\alpha\in Q\mid\forall f\in W,\alpha(f)=0\}.
\end{displaymath}
In particular, if $\ell=1$ and $W=\langle f\rangle$ then one writes $\Ann(W)=\Ann(f)$.
\end{definition}
\noindent Let $A=Q/\Ann(f)$, where $f$ is homogeneous. By construction $A$ is a graded Artinian $\K$-algebra; moreover $A$ is also Gorenstein.
\begin{theorem}[{\cite[\S60ff]{M:FHS}}, {\cite[theorem 2.1]{M:W}}]\label{th1.1}
Let $I$ be a homogeneous ideal of $Q$ such that $A=Q/I$ is a graded Artinian $\K$-algebra. Then $A$ is Gorenstein of socle degree $d$ if and only if there exist $d\geq1$ and $f\in R_d$ such that $A\cong Q/\Ann(f)$.
\end{theorem}
\begin{remark}
Using the above notations, if $\Ann(f)_1=0$ then $A$ is called the \emph{standard graded Artinian Gorenstein} (\emph{SGAG}, for short) \emph{$\K$-algebra associated to $f$}. The socle degree $d$ of $A$ is the degree of $f$ and the codimension is $m$.
\end{remark}

\subsection{Full Perazzo algebras}

\begin{definition}
Let $\displaystyle A=\bigoplus_{i=0}^dA_i$ be a SGAG $\K$-algebra. $A$ is a \emph{standard bigraded Artinian Gorenstein algebra} (\emph{SBAG algebra}, for short) if it is \emph{bigraded}, \emph{i.e.}: 
\begin{displaymath}
A_d=A_{(d_1,d_2)}\cong\K,\,A_i=\bigoplus_{h=0}^i A_{(i,h-i)}\,\text{\,for each\,}\,i\in\{0,\dotsc,d-1\}.
\end{displaymath}
Since $A$ is a Gorenstein ring, the pair $(d_1,d_2)$ is said the \emph{socle bidegree of $A$}.
\end{definition}
\begin{definition}
A homogeneous ideal $I$ of $S$ is a \emph{bihomogeneous ideal} if:
\begin{displaymath}
I=\bigoplus_{i,j=0}^{+\infty}I_{(i,j)},\text{\,where\,}\,\forall i,j\in\N_{\geq0},\,I_{(i,j)}=I\cap S_{(i,j)}.
\end{displaymath}
\end{definition}
\noindent Let $f\in S_{(d_1,d_2)}$, then $\Ann(f)$ is a bihomogeneous ideal and using Theorem \ref{th1.1}, $A=T/\Ann(f)$ is a SBAG $\K$-algebra of socle bidegree $(d_1,d_2)$ (and codimension $m+n$).
\begin{definition}
A bihomogeneous polynomial
\begin{displaymath}
f=\sum_{i=1}^nx_ig_i\in\K\left[x_1,\dotsc,x_m,u_1,\dotsc,u_n\right]_{(1,d-1)},
\end{displaymath}
where $g_i$'s are are algebraically dependent but linearly independent monomials of degree $d-1$, is called a \emph{Perazzo polynomial} of degree $d$ (see \cite[Definition 3.1]{CGIZ}).\smallskip

\noindent If the monomials $g_i$ form a basis of $\K\left[u_1,\dotsc,u_n\right]_{d-1}$ then $f$ is a \emph{full Perazzo polynomial} of degree $d$ (see \cite[Definition 3.2]{CGIZ}).\smallskip

\noindent The associated algebra $A=T/\Ann(f)$ is called (\emph{full}) \emph{Perazzo algebra}.
\end{definition}
\begin{remark}
One needs $\displaystyle m\leq\binom{n+d-2}{d-1}$ otherwise the $g_i$'s cannot be linearly independent.\smallskip

\noindent From now on, I assume that $m$ satisfies this condition.
\end{remark}

\subsection{Lefschetz properties}

\noindent Let $\displaystyle A=\bigoplus_{i=0}^d A_i$ be a graded Artinian $\K$-algebra.
\begin{definition}
If there exists an $L\in A_1$ such that:
\begin{enumerate}[a)]
\item The multiplication map $\cdot L\colon A_i\to A_{i+1}$ is of maximal rank for all $i$, then $A$ has the \emph{Weak Lefschetz Property} (\emph{WLP}, for short);
\item The multiplication map $\cdot L^k\colon A_i\to A_{i+k}$ is of maximal rank for all $i$ and $k$, then $A$ has the \emph{Strong Lefschetz Property} (\emph{SLP}, for short);
\end{enumerate}
\end{definition}
\noindent About Perazzo algebras with minimal Hilbert vector, the following theorem holds.
\begin{theorem}\label{th1.2}
Let $m\geq n\in\N_{\geq 2}$, $d\in\N_{\geq1}$, let $f$ be a Perazzo polynomial and let $A=T/\Ann(f)$ be the associated SGAG algebra. Suppose that $A$ has minimal Hilbert vector then: 
\begin{enumerate}[a)]
\item let $n=2$. $A$ has WLP if and only if $d\geq2(m-1)$ (\cite[Corollary 4.11]{MR:P}).
\item let $n\in\N_{\geq3}$. If the Hilbert vector of $A$ is unimodal and $m\leq\binom{d+n-s-1}{n-1}$ then $A$ has WLP, where $s=\displaystyle\left\lfloor\frac{d}{2}\right\rfloor$ (\cite[Theorem 4.5]{M:MR}).
\end{enumerate}
\end{theorem}

\section{On the locus of socle degree $d$ SGAG algebras of codimension $n$}
\markboth{On the locus of socle degree $d$ SGAG algebras of codimension $n$}{On the locus of socle degree $d$ SGAG algebras of codimension $n$}

\subsection{Preliminaries}

\noindent In order to state and prove some lemma, I prove the following known propositions for which I am not able to indicate a reference.
\begin{proposition}\label{prop2.1}
Let $n\in\N_{\geq1}$, let $\{I,J\}$ be a partition of $\{0,\dotsc,n\}$. The locally closed set
\begin{displaymath}
V_{I,J}=\left\{\left[a_0:\dotsc:a_n\right]\in\bP^n\mid\forall i\in I,j\in J,a_i=0,a_j\neq0\right\}
\end{displaymath}
is a smooth, affine and irreducible variety of dimension $\# J-1$, with respect to Zariski topology.
\end{proposition}
\begin{prf}
For $n=1$ one has a point hence there is nothing to prove. Let $n\geq2$, up to a change of coordinates, one may assume $J=\{0,\dotsc,k\}$ and $I=\{k+1,\dotsc,n\}$, \emph{i.e.} one may assume that $V_{I,J}$ is the set
\begin{displaymath}
\left\{\left[a_0:\dotsc:a_k:0:\dotsc:0\right]\in\bP^n\mid a_0\cdot\dotsc\cdot a_k\neq0\right\}
\end{displaymath}
which is isomorphic to
\begin{displaymath}
\left\{\left[a_0:\dotsc:a_k\right]\in\bP^k\mid a_0\cdot\dotsc\cdot a_k\neq0\right\}.
\end{displaymath}
Let
\begin{displaymath}
U_0=\left\{\left[a_0:\dotsc:a_k\right]\in\bP^k\mid a_0\neq0\right\},
\end{displaymath}
it is known $U_0\cong\A^k$; since $V_{I,J}\subseteq U_0$ then it is isomorphic to
\begin{displaymath}
V=\left\{\left(a_1,\dotsc,a_k\right)\in\A^k\mid a_1\cdot\dotsc\cdot a_k\neq0\right\}.
\end{displaymath}
Using the \emph{Rabinowitsch's trick}, $V$ is isomorphic to the affine irreducible variety
\begin{displaymath}
\left\{\left(a_1,\dotsc,a_k,t\right)\in\A^{k+1}\mid a_1\cdot\dotsc\cdot a_k\cdot t=1\right\}
\end{displaymath}
which is smooth by \emph{Jacobian criterion} (\cite[Exercises 13.2.E and 13.2.J]{FOAG}) and has dimension $k$.
\end{prf}
\begin{proposition}\label{prop2.4}
Let $f\colon\V\to\W$ be a linear map of vector spaces. Let
\begin{displaymath}
\widehat{f}\colon\left[\uv\right]\in\bP(\V)\dashrightarrow\left[f\left(\uv\right)\right]\in\bP(\Im(f))
\end{displaymath}
be the surjective rational map induced by $f$. Then $\widetilde{f}=\widehat{f}_{\vert\bP(\V)\setminus\bP(\ker(f))}$ is a surjective smooth morphism of relative dimension $\dim\ker(f)-1$.
\end{proposition}
\begin{prf}
By the previous proposition, $\bP(\V)\setminus\bP(\ker(f))$ is a smooth irreducible affine variety. Thus, the differential $d\widehat{f}$ of $\widehat{f}$ restricted to this variety is given by Jacobian matrix of $f$ (\cite[Example 6.4]{G:W}), in other words $d\widehat{f}$ is a surjective linear map at any point of $\bP(\V)\setminus\bP(\ker(f))$, this means that $\widehat{f}$ is smooth of relative dimension $\dim\ker(f)-1$ by \cite[Proposition III.10.4]{H:RC}.
\end{prf}
\noindent In this section, from now on let $m\in\N_{\geq2}$ and let $d\in\N_{\geq1}$, unless otherwise indicated.
\begin{remark}
I carry on the set $S(m,d)=\{0,\dotsc,d\}^m$ the \emph{lexicographic order} $\preceqlex^m$.\smallskip

\noindent I recall that $\left(a_1,\dotsc,a_m\right)\preceqlex^m\left(b_1,\dotsc,b_m\right)$ if and only if either $a_1\preceqlex^1b_1$ or $a_1=b_1$ and $\left(a_2,\dotsc,a_m\right)\preceqlex^{m-1}\left(b_2,\dotsc,b_m\right)$, where $\preceqlex^1$ is the natural order of $\N_{\geq0}$.
\end{remark}
\noindent Let $\ui=\left(i_1,\dotsc,i_m\right)\in S(m,d)$, $\left|\ui\right|=i_1+\dotsc+i_m$, let
\begin{equation}\label{eq2}
T(m,d)=\left\{\ui\in S(m,d):\left|\ui\right|=d\right\};
\end{equation}
on this set I consider the order induced by $\left(S(m,d),\preceqlex^m\right)$, unless otherwise indicated. I shall use the following functions in the next:
\begin{gather}
\forall k\in\{1,\dotsc,m\},\,\pr_k^{(m,d)}\colon\ui\in S(m,d)\to i_k\in\{0,\dotsc,d\}\label{eq4};\\
D_k^{(m,d)}\colon\ui\in S(m,d)\to S(m,d)\cup\{-\infty\}\ni\begin{cases}
\left(i_1,\dots,i_{k-1},i_k-1,i_{k+1},\dotsc,i_m\right)\iff\pr_k^{(m,d)}\left(\ui\right)\geq1\\
-\infty\iff\pr_k^{(m,d)}\left(\ui\right)=0
\end{cases}, \label{eq6}\\
d_k^{(m,d)}\colon\left(i_1,\dots,i_m\right)\in S(m,d)\to\left(i_1,\dots,i_{k-1},i_{k+1},\dotsc,i_m\right)\in S(m-1,d) \notag.
\end{gather}
\begin{remark}
With abuse of the notations, I shall skip any reference to $m$ and $d$ where I write the previous functions.
\end{remark}

\subsection{The standard locus}

\noindent For the aims of this paper, I consider the following set
\begin{equation}\label{eq1}
U^{(m,d)}_{st}=\left\{[f]\in\bP\left(R_d\right)\mid\forall\lambda\in Q_1,\,\lambda(f)\neq0\right\}.
\end{equation}
\noindent Since $[f]\in U^{(m,d)}_{st}$ if and only if $A/\Ann(f)$ is a socle degree $d$ SGAG $\K$-algebra of codimension $m$, I call this set \emph{standard locus of socle degree $d$ and codimension $m$}.
\begin{lemma}\label{lem2.1}
$U^{(m,d)}_{st}$ is union of affine locally closed subsets of $\bP\left(R_d\right)$ with respect to Zariski topology.
\end{lemma}
\begin{prf}
Recall that $\displaystyle\dim_{\K}R_d=\binom{m+d-1}{d}$ which is denoted by $\tau$, for simplicity. One puts
\begin{equation}\label{eq3}
x^{\ui}=x_1^{i_1}\cdot\dotsc\cdot x_m^{i_m}\,\text{\,and\,}\,x^{-\infty}=0.
\end{equation}
By assumptions
\begin{displaymath}
\forall f\in R_d,\,\lambda\in Q_1,\,\dot\exists\left(a_{\ui}\right)\in\K^{\tau}\setminus\left\{\uzero^{\tau}\right\},\left(b_j\right)\in\K^m\setminus\left\{\uzero^m\right\}\mid f=\sum_{\ui\in T(m,d)}a_{\ui}x^{\ui},\,\lambda=\sum_{j=1}^mb_jX_j;
\end{displaymath}
in particular one may confuse $f$ and $\lambda$ with their ordered coefficients. By definition
\begin{displaymath}
\lambda(f)=\sum_{j=1}^m\sum_{\ui\in T(m,d)}a_{\ui}b_jx^{D_j\left(\ui\right)}\in R_{d-1}
\end{displaymath}
hence $[f]\in U^{(m,d)}_{st}$ if and only if
\begin{displaymath}
\forall\left(b_j\right)\in\K^n\setminus\left\{\uzero^n\right\},\,\exists\ui^{\prime}\in T(m,d-1)\mid\sum_{j=1}^m\sum_{\ui\in T(m,d)\mid D_j\left(\ui\right)=\ui^{\prime}}a_{\ui}b_j\neq0.
\end{displaymath}
In other words, fixed $[f]\in U^{(m,d)}_{st}$:
\begin{enumerate}
\item[st.A)] for each $j\in\{1,\dotsc,m\}$ there exists $\ui\in T(m,d)$ such that $a_{\ui}\neq0$ and $D_j\left(\uj\right)\neq-\infty$;
\item[st.B)] for each $\ui^{\prime}\in T(m,d-1)$ either there exist unique $\ui\in T(m,d)$ and ${j\in\{1,\dotsc,m\}}$ such that $a_{\ui}\neq0$ and $D_j\left(\ui\right)=\ui^{\prime}$ or they do not exist;
\item[st.C)] there do not exist $\ui_1\neq \ui_2\in T(m,d)$ and $j_1\neq j_2\in\{1,\dotsc,m\}$ such that ${a_{\ui_1},a_{\ui_2}\neq0}$ and $D_{j_1}\left(\ui_1\right)=D_{j_2}\left(\ui_2\right)$.
\end{enumerate}
Thus $U^{(m,d)}_{st}$ is the union of sets $U^{(m,d)}_{\wl}$ of the polynomials $\displaystyle f=\sum_{\ui\in T(m,d)}a_{\ui}x^{\ui}\in R_d$ such that the coefficients of $f$ satisfy conditions (st.A), (st.B) and (st.C) and each $D_j\left(\ui\right)$, with $a_{\ui}\neq0$, belongs to $\wl$; where $\wl$ is a fixed element of the following set:
\begin{gather*}
T(m,d)_{st}=\left\{\wl=\left\{\ui^{\prime}_1,\dotsc,\ui^{\prime}_p\right\}\in\sP\left(T(m,d-1)\right)\mid\right.\\
\left.\forall k\in\{1,\dotsc,p\},\,\dot\exists j_k\in\{1,\dotsc,m\},\ui_k\in T(m,d):D_{j_k}\left(\ui_k\right)=\ui_k^{\prime}\right.\\
\left.\forall j\in\{1,\dotsc,m\},\,\exists k_j\in\{1,\dotsc,p\},\ui_j\in T(m,d):D_j\left(\ui_j\right)=\ui^{\prime}_{k_j}\right\}.
\end{gather*}
By Proposition \ref{prop2.1}, each set $U^{(m,d)}_{\wl}$ is smooth, irreducible, affine and locally closed. Thus the claim holds.
\end{prf}
\begin{corollary}
Let $\displaystyle\left[\sum_{\ui\in T(m,d)}a_{\ui}x^{\ui}\right]\in U^{(m,d)}_{st}$. Then for any $\ui_1\neq\ui_2\in T(m,d)$ such that $a_{\ui_1},a_{\ui_2}\neq0$ it results $\deg\left(\gcd\left(x^{\ui_1},x^{\ui_2}\right)\right)\leq d-2$.
\end{corollary}
\begin{lemma}\label{lem2.3}
Without changing the notations from the proof of Lemma \ref{lem2.1}, the affine smooth sets $U^{(m,d)}_{\wl}$ are the irreducible components of $U^{(m,d)}_{st}$ of dimensions $\#\wl-1$.
\end{lemma}
\begin{prf}
By Proposition \ref{prop2.1}, one needs only to prove that $U^{(m,d)}_{\wl}$'s are the irreducible components of $U^{(m,d)}_{st}$. Without changing the notations from Lemma \ref{lem2.1}:
\begin{displaymath}
\forall\wl_1\neq\wl_2\in T(m,d)_{st},\,U^{(m,d)}_{\wl_1}\nsubseteq U^{(m,d)}_{\wl_2},U^{(m,d)}_{\wl_2}\nsubseteq U^{(m,d)}_{\wl_1},U^{(m,d)}_{\wl_1}\cap U^{(m,d)}_{\wl_2}=\emptyset
\end{displaymath}
and all this proves the claim.
\end{prf}
\noindent Moreover, by Proposition \ref{prop2.1}, the following corollary hold.
\begin{corollary}\label{cor2.1}
$U^{(m,d)}_{st}$ is a smooth subset of $\bP\left(R_d\right)$.
\end{corollary}
\begin{corollary}
For any $m\in\N_{\geq2}$ and $d\in\N_{\geq1}$
\begin{displaymath}
\left\{\left[a_1x_1^d+\dotsc+a_mx_m^d\right]\in\bP\left(R_d\right)\mid a_1\cdot\dotsc\cdot a_m\neq0\right\}\subseteqq U^{(m,d)}_{st}
\end{displaymath}
is an irreducible component of $U^{(m,d)}_{st}$ of dimension $m-1$; in particular $U^{(m,d)}_{st}$ is not empty.
\end{corollary}
\noindent From now on, let $M_1,\dotsc,M_{\tau(n,d-1)}$ be the monic monomial basis of $\K\left[u_1,\dotsc,u_n\right]_{d-1}$ where $\displaystyle{\tau(n,d)=\binom{n+d-1}{d}}$. For sake of simplicity, one puts $\left(\ui,\uj\right)=\left(i_1,\dotsc,i_{\tau(n,d-1)},j_1,\dotsc,j_n\right)$ and $x^{\ui}u^{\uj}$ or $(xu)^{\left(\ui,\uj\right)}$ both mean $x^{i_1}\cdot\dotsc\cdot x_{\tau(n,d-1)}^{i_{\tau(n,d-1)}}u_1^{j_1}\cdot\dotsc\cdot u_n^{j_n}$. Let $Z_{(n,d),1}$ be the following set
\begin{displaymath}
\left\{k\in\{1,\dotsc,\tau(n,d-1)\}\mid u_n/x_kM_k\right\}=\left\{k_1<\dotsc<k_t\right\},
\end{displaymath}
where $u_n/x_kM_k$ means $u_n$ divides $x_kM_k$. Let
\begin{displaymath}
\varphi_{(n,d)}\colon\K\left[x_1,\dotsc,x_{\tau(n,d-1)},u_1,\dotsc,u_n\right]_d\to\K\left[x_1,\dotsc,x_{\tau(n-1,d-1)},u_1,\dotsc,u_{n-1}\right]_d
\end{displaymath}
be the linear morphism of vector spaces such that
\begin{itemize}
\item $x^{\ui}u^{\uj}\in\ker\left(\varphi_{(n,d)}\right)$ if and only if $\pr_n\left(\uj\right)\geq1$ or there exists $h\in Z_{(n,d),1}$ such that $x_h/x^{\ui}$;
\item
$\varphi_{(n,d)}\left(x^{\ui}u^{\uj}\right)=x^{\ui^{\prime}}u^{\uj\prime}$, where $x^{\ui}u^{\uj}\notin\ker(\varphi_{(n,d)})$, ${\left(d_{k_1}\circ\dotsc\circ d_{k_t}\right)\left(\ui\right)=\ui^{\prime}}$ and $d_n\left(\uj\right)=\uj^{\prime}$.
\end{itemize}
\noindent By construction $\varphi_{(n,d)}$ is surjective hence
\begin{displaymath}
\dim\ker\left(\varphi_{(n,d)}\right)=\tau(n,d-1)-\tau(n-1,d-1)+1.
\end{displaymath}
Let
\begin{displaymath}
\whvarphi_{(n,d)}\colon\bP\left(\K\left[x_1,\dotsc,x_{\tau(n,d-1)},u_1,\dotsc,u_n\right]_d\right)\dashrightarrow\bP\left(\K\left[x_1,\dotsc,x_{\tau(n-1,d-1)},u_1,\dotsc,u_{n-1}\right]_d\right)
\end{displaymath}
be the surjective rational map induced by $\varphi_{(n,d)}$.
\begin{proposition}\label{prop2.2}
The restriction $\wtvarphi_{(n,d)}$ of $\whvarphi_{(n,d)}$ to $U^{(n+\tau(n,d-1),d)}_{st}$ is a smooth surjective morphism onto $U^{(n-1+\tau(n-1,d-1),d)}_{st}$ of relative dimension ${\dim\ker\left(\varphi_{(n,d)}\right)-1}$.
\end{proposition}
\begin{prf}
By the proof of Lemma \ref{lem2.1}, $\bP\left(\ker\left(\whvarphi_{(n,d)}\right)\right)\cap U^{(n+\tau(n,d-1),d)}_{st}=\emptyset$ \emph{i.e.} $\wtvarphi_{(n,d)}$ is a regular morphism. Let $[f]\in U^{(n-1+\tau(n-1,d-1),d)}_{st}$, it results
\begin{displaymath}
\left[f+x_{k_1}^d+\dotsc+x_{k_t}^d+u_n^d\right]\in U^{(n+\tau(n,d-1),d)}_{st}\,\text{\,and\,}\,\wtvarphi_{(n,d)}\left(\left[f+x_{k_1}^d+\dotsc+x_{k_t}^d+u_n^d\right]\right)=[f]
\end{displaymath}
\emph{i.e.} $\wtvarphi_{(n,d)}$ is surjective onto $U^{(n-1+\tau(n-1,d-1),d)}_{st}$. Since the irreducible components of $U^{(n+\tau(n,d-1),d)}_{st}$ are smooth and affine (Proposition \ref{prop2.1}), by Proposition \ref{prop2.4} one has the claim.
\end{prf}
\noindent Let
\begin{gather*}
h_{(n,d)}\colon\ul\in T\left(n,d\right)\to\max\left\{k\in\{1,\dotsc,n\}\mid\pr_k\left(\ul\right)\geq1\right\}\in\N_{\geq0},\\
c_{(n,d)}\colon\ul\in T\left(n,d\right)\to\ul^{\prime}\in T\left(n,d-1\right)
\end{gather*}
where $\ell^{\prime}$ is defined as it follows:
\begin{displaymath}
\pr_{h_{(n,d-1)}\left(\ul^{\prime}\right)}\left(\ul^{\prime}\right)=\pr_{h_{(n,d)}\left(\ul\right)}\left(\ul\right)-1,\,d_{h_{(n,d-1)}\left(\ul^{\prime}\right)}\left(\ul^{\prime}\right)=d_{h_{(n,d)}\left(\ul\right)}\left(\ul\right),
\end{displaymath}
and let
\begin{gather*}
c_{(n,d),0}\colon\ul^{\prime}\in T\left(n,d-1\right)\to\min c^{-1}_{(n,d)}\left(\ul^{\prime}\right)\in T\left(n,d\right),\\
Z_{(n,d),2}=\left\{k\in\left\{1,\dotsc,\tau(n,d-1)\right\}\mid M_k=u^{\ul}\,\text{\,for some\,}\,\ul\in\Im\left(c_{(n,d-1),0}\right)\right\}.
\end{gather*}
\noindent Let
\begin{displaymath}
\psi_{(n,d)}\colon\K\left[x_1,\dotsc,x_{\tau(n,d-1)},u_1,\dotsc,u_n\right]_d\to\K\left[x_1,\dotsc,x_{\tau(n,d-2)},u_1,\dotsc,u_n\right]_{d-1}
\end{displaymath}
be the linear morphism of vector spaces such that
\begin{itemize}
\item $x^{\ui}u^{\uj}\in\ker\left(\psi_{(n,d)}\right)$ if and only if at least one of the following conditions occurs:
\begin{enumerate}[I)]
\item there exists $k\in\{1,\dotsc,\tau(n,d-1)\}\setminus Z_{(n,d),2}$ such that $x^k/x^{\ui}$,
\item $\left(\ui,\uj\right)\notin\Im\left(c_{(n+\tau(n,d-1),d),0}\right)$,
\item $\uj\notin\Im\left(c_{(n,d-1),0}\right)$;
\end{enumerate}
\item $\psi_{(n,d)}\left(x^{\ui}u^{\uj}\right)=(xu)^{c_{(n+\tau(n,d-1),d)}\left(\ui,\uj\right)}$ if $x^{\ui}u^{\uj}\notin\ker\left(\psi_{(n,d)}\right)$.
\end{itemize}
\noindent By construction $\psi_{(n,d)}$ is surjective hence
\begin{displaymath}
\dim\ker\left(\psi_{(n,d)}\right)=\tau(n,d-1)-\tau(n,d-2).
\end{displaymath}
Let
\begin{displaymath}
\whpsi_{(n,d)}\colon\bP\left(\K\left[x_1,\dotsc,x_{\tau(n,d-1)},u_1,\dotsc,u_n\right]_d\right)\dashrightarrow\bP\left(\K\left[x_1,\dotsc,x_{\tau(n,d-2)},u_1,\dotsc,u_n\right]_{d-1}\right)
\end{displaymath}
be the surjective rational map induced by $\psi_{(n,d)}$.\smallskip

\noindent Repeating the proof of Proposition \ref{prop2.2} one has the following result.
\begin{proposition}\label{prop2.3}
The restriction $\wtpsi_{(n,d)}$ of $\whpsi_{(n,d)}$ to $U^{(n+\tau(n,d-1),d)}_{st}$ is a smooth surjective morphism onto $U^{(n+\tau(n,d-2),d)}_{st}$ of relative dimension ${\dim\ker\left(\psi_{(n,d)}\right)-1}$.
\end{proposition}

\subsection{The Full Perazzo locus}

\noindent Let $U^{(n,d)}_{FP}$ be the subset of $U^{(n+\tau(n,d-1),d)}_{st}$ determined by the union of the following sets
\begin{gather}
\left\{\left[a_1x_1M_1+\dotsc+a_{\tau(n,d-1)}x_{\tau(n,d-1)}M_{\tau(n,d-1)}\right]\in\bP\left(\K\left[x_1,\dotsc,x_{\tau(n,d-1)},u_1,\dotsc,u_n\right]_{(1,d-1)}\right)\mid\right. \notag\\
\left.a_1\cdot\dotsc\cdot a_{\tau(n,d-1)}\neq0\right\} \label{eq5}
\end{gather}
up to a change of coordinates. I call this set \emph{full Perazzo locus}.
\begin{lemma}\label{lem2.2}
For each $n,d\in\N_{\geq2}$, $U^{(n,d)}_{FP}$ has pure dimension $\tau(n,d-1)-1$. And it does not exist an irreducible component of $U^{(n+\tau(n,d-1),d)}_{st}$ whose dimension is less than $\tau(n,d-1)-1$.
\end{lemma}
\noindent In the following proof I use the smooth surjective morphisms $\wtvarphi_{(n,d)}$ and $\wtpsi_{(n,d)}$ defined above.
\begin{prf}
By Proposition \ref{prop2.1}, $U^{(n,d)}_{FP}$ has pure dimension $\tau(n,d-1)-1$.\smallskip

\noindent Let $(n,d)=(2,2)$ then $\dim U^{(2,2)}_{FP}=1$ and there is nothing to prove, since it does not exist a $0$-dimensional irreducible component of $U^{(2+\tau(2,1),2)}_{st}$ by Proposition \ref{prop2.1} and proof of Lemma \ref{lem2.1}. Let $(n,d)=(3,2)$, by construction there exists an irreducible component of $U^{(3,2)}_{FP}$ which is transformed by $\wtvarphi_{(3,2)}$ in an irreducible component of $U^{(2,2)}_{FP}$ \emph{i.e.} $\dim\wtvarphi_{(3,2)}\left(U^{(3,2)}_{FP}\right)=\dim U^{(2,2)}_{FP}$. By Proposition \ref{prop2.2}, there is no one irreducible component of $U^{(3+\tau(3,1),2)}_{st}$ which is contracted by $\wtvarphi_{(3,2)}$ to a point. From all this, by \cite[Definition at page 268]{H:RC}
\begin{gather*}
\forall V\,\text{\,irreducible component of\,}\,U^{(3+\tau(3,1),2)}_{st},\\
\dim V\geq\dim U^{(2,2)}_{FP}+\dim\ker\left(\wtvarphi_{(3,2)}\right)-1=\tau(3,1)-1;
\end{gather*}
in particular, the equality holds for the irreducible components of $U^{(3,2)}_{FP}$; \emph{i.e.} the statement holds for $(n,d)=(3,2)$. By inductive hypothesis, let the statement holds for $(n,2)$ where $n\in\N_{\geq3}$, by construction, there exists an irreducible component of $U^{(n+1,2)}_{FP}$ which is transformed by $\wtvarphi_{(n+1,2)}$ in an irreducible component of $U^{(n,2)}_{FP}$. Repeating the previous reasoning one has
\begin{gather*}
\forall V\,\text{\,irreducible component of\,}\,U^{(n+1+\tau(n+1,1),2)}_{st},\\
\dim V\geq\dim U^{(n,2)}_{FP}+\dim\ker\left(\wtvarphi_{(n+1,2)}\right)-1=\tau(n+1,1)-1
\end{gather*}
and the equality holds for the irreducible components of $U^{(n,2)}_{FP}$. By \emph{Mathematical Induction}, the statement holds for each pair $(n,2)$. On another hand, fix $n\in\N_{\geq2}$ let $d=3$, by construction there exists an irreducible component of $U^{(n,3)}_{FP}$ which is transformed by $\wtpsi_{(n,3)}$ in an irreducible component of $U^{(n,2)}_{FP}$ \emph{i.e.} $\dim\wtpsi_{(n,3)}\left(U^{(n,3)}_{FP}\right)=U^{(n,2)}_{FP}$. By Proposition \ref{prop2.3} and the previous reasoning, there is no one irreducible component of $U^{(n+\tau(n,2),3)}_{st}$ which is contracted by $\wtpsi_{(n,3)}$ to an irreducible subset of dimension less than $\dim U^{(n,2)}_{FP}$. From all this, by \cite[Definition at page 268]{H:RC}:
\begin{gather*}
\forall V\,\text{\,irreducible component of\,}\,U^{(n+\tau(n,2),3)}_{st},\\
\dim V\geq\dim U^{(n,3)}_{FP}+\dim\ker\left(\wtpsi_{(n,3)}\right)-1=\tau(n,2)-1,
\end{gather*}
in particular, the equality holds for the irreducible components of $U^{(n,3)}_{FP}$; \emph{i.e.} the statement holds for each $n\in\N_{\geq2}$ and $d=3$. By inductive hypothesis, let the statement holds for $(n,d)$ where $d\in\N_{\geq3}$, by construction, there exists an irreducible component of $U^{(n,d+1)}_{FP}$ which is transformed by $\wtvarphi_{(n,d+1)}$ in an irreducible component of $U^{(n,d)}_{FP}$. Repeating the previous reasoning one has
\begin{gather*}
\forall V\,\text{\,irreducible component of\,}\,U^{(n+\tau(n,d),d+1)}_{st},\\
\dim V\geq\dim U^{(n,d)}_{FP}+\dim\left(\ker\wtpsi_{(n,d)}\right)-1=\tau(n,d-1)-1
\end{gather*}
and by \emph{Mathematical Induction}, the statement holds for each pair $(n,d)$.
\end{prf}

\section{CW-complexes and simplicial complexes}
\markboth{CW-complexes and simplicial complexes}{CW-complexes and simplicial complexes}

\subsection{Abstract finite simplicial complexes}

\begin{definition}
Let $V=\left\{x_1,\dotsc,x_m\right\}$ be a finite set. An \emph{abstract finite simplicial complex $\Delta$ with vertex set $V$} is a subset of $\sP(V)$ such that
\begin{enumerate}[1)]
\item $\forall x\in V\Rightarrow\{x\}\in\Delta$,
\item $\forall\sigma\in\Delta,\tau\subsetneqq\sigma,\tau\neq\emptyset\Rightarrow\tau\in\Delta$.
\end{enumerate}
\end{definition}
\noindent The elements $\sigma$ of $\Delta$ are called \emph{faces} or \emph{simplices}; a face with $q+1$ vertices is called \emph{$q$-face} or \emph{face of dimension $q$} and one writes $\dim\sigma=q$; the maximal faces (with respect to the inclusion) are called \emph{facets}; if all facets have the same dimension $d\geq1$ then one says $\Delta$ has \emph{pure dimension $d$}. The set $\Delta^k$ of faces of dimension at most $k$ is called \emph{$k$-skeleton of $\Delta$}; instead, the set of faces of $\Delta$ of dimension $k$ is denoted as $F^k(\Delta)$. $\sP(V)$ is called \emph{simplex} (of dimension $m-1$).
\begin{remark}[{cfr. \cite[Remark 3.4]{CGIM}}]\label{rem3.1}
There is a natural bijection, described in \cite[Section 5.1]{B:H}, between the square-free monomials, of degree $d$, in the variables $x_1,\dotsc,x_m$ and the $(d-1)$-faces of the simplex $\sP(V)$, with vertex set $V=\left\{x_1,\dotsc,x_m\right\}$. In fact, a square-free monomial $g=x_{i_1}\cdot\dotsc\cdot x_{i_d}$ corresponds to the subset $\left\{x_{i_1},\dotsc,x_{i_d}\right\}$ of $V$. \emph{Vice versa}, to any subset $W$ of $V$ with $d$ elements one associates the free square monomial $\displaystyle m_W=\prod_{x_i\in W}x_i$ of degree $d$.
\end{remark}

\subsection{CW-complexes}

For the topological background, I refer to \cite{H:AE}. I start by fixing some notations. 
\begin{definition}
Let $k\in\N_{\geq1}$. A topological space $e^k$ homeomorphic to the open (unitary) ball
\begin{displaymath}
\B^k=\left\{\left(x_1,\dotsc,x_k\right)\in\R^k\mid x_1^2+\cdot\dotsc\cdot+x_k^2<1\right\}
\end{displaymath}
of dimension $k$ (with the natural topology induced by $\R^{k+1}$) is called a \emph{$k$-cell}. Its boundary $\partial e^k$, \emph{i.e.} the \emph{$(k-1)$-dimensional sphere}, will be denoted by
\begin{displaymath}
\bS^{k-1}=\left\{\left(x_1,\dotsc,x_k\right)\in\R^k\mid x_1^2+\cdot\dotsc\cdot+x_k^2=1\right\}
\end{displaymath}
and its closure $\overline{e}^k$, \emph{i.e.} the closed (unitary) $k$-dimensional disk, will be denoted by
\begin{displaymath}
\D^k=\left\{\left(x_1,\dotsc,x_k\right)\in\R^k\mid x_1^2+\cdot\dotsc\cdot+x_k^2\leq1\right\}.
\end{displaymath} 
\end{definition}
\noindent I recall the following 
\begin{definition}
A \emph{CW-complex} is a topological space $X$ constructed in the following way:
\begin{enumerate}[a)]
\item There exists a fixed and discrete set of points $X^0\subseteq X$, whose elements are called \emph{$0$-cells};
\item Inductively, the \emph{$k$-skeleton $X^k$ of $X$} is constructed from $X^{k-1}$ by attaching $k$-cells $e_{\alpha}^k$ (with index set $A_k$) via continuous maps $\varphi^k_{\alpha}\colon\bS^{k-1}_{\alpha}\to X^{k-1}$ (the \emph{attaching maps}). This means $X^k$ is a quotient of $\displaystyle Y^k=\bigcup_{\alpha\in A_k}\D^k_{\alpha}\cup X^{k-1}$ under the identification $x\sim\varphi_{\alpha}(x)$ for each $x\in\bS^{k-1}_{\alpha}$, and $e^k_{\alpha}$ is the homeomorphic image of $\B^k_{\alpha}$ under the quotient map;
\item $\displaystyle X=\bigcup_{k\in\N_{\geq0}}X^k$.
\end{enumerate}
\noindent The \emph{elements of the $k$-skeleton} are the (closure of the) attached $k$-cells. A subset $C$ of $X$ is closed if and only if $C\cap X^k$ is closed for any $k\in\N_{\geq0}$ (\emph{closed weak topology}).
\end{definition}
\begin{definition}
A subset $Z$ of a CW-complex $X$ is a \emph{CW-subcomplex} if it is the union of cells of $X$, such that the closure of each cell is in $Z$.
\end{definition}
\begin{definition}
A CW-complex is \emph{finite} if it consists of a finite number of cells.
\end{definition}
\noindent I shall be interested mainly in finite CW-complexes.
\begin{example}[\emph{Geometric realization of an abstract simplicial complex}]
It is an obvious fact that to any finite simplicial complex $\Delta$ one can associate a finite CW-complex $\widetilde{\Delta}$ via the \emph{geometric realization} of $\Delta$ as a simplicial complex (as a topological space) $\widetilde{\Delta}$.
\end{example}
\noindent In what follows I shall always identify abstract simplicial complexes with their corresponding simplicial complexes.

\section{CW-complexes and homogeneous polynomials}
\markboth{CW-complexes and homogeneous polynomials}{CW-complexes and homogeneous polynomials}

\begin{construction}[{see \cite[Chapter 5]{B:H}}]
In Remark \ref{rem3.1}, one saw that to any degree $d$ square-free monomial $x_{i_1}\cdot\dotsc\cdot x_{i_d}\in R_d$ one can associate the $(d-1)$-face $\left\{x_{i_1},\dotsc,x_{i_d}\right\}$ of the abstract $(m-1)$-dimensional simplex $\Delta(m)=\sP\left(\left\{x_1,\dotsc,x_m\right\}\right)$, and \emph{vice versa}: let
\begin{displaymath}
\rho_d=\left\{f\in R_d\mid f\neq0\,\text{\,is a square-free monic monomial}\right\},
\end{displaymath}
one has a bijection
\begin{align*}
\sigma_d\colon\rho_d &\to F^d(\Delta(m))\\
x_{i_1}\cdot\dotsc\cdot x_{i_d} & \mapsto\left\{x_{i_1},\dotsc,x_{i_d}\right\}.
\end{align*}
Alternatively, one can associate to $x_{i_1}\cdot\dotsc\cdot x_{i_d}$ the element of the $(d-1)$-skeleton $\overline{\left\{x_{i_1},\dotsc,x_{i_d}\right\}}\in\widetilde{\Delta(m)}^{d-1}$, so one has a bijection
\begin{align*}
\sigma_d\colon\rho_d & \to\widetilde{\Delta(m)}^{d-1}\setminus\widetilde{\Delta(m)}^{d-2}\\
x_{i_1}\cdot\dotsc\cdot x_{i_d} & \mapsto\overline{\left\{x_{i_1},\dotsc,x_{i_d}\right\}}
\end{align*}
between the square-free monomials and the $(d-1)$-faces of the (topological) simplex $\widetilde{\Delta(m)}$.
\end{construction}
\begin{construction}[{cfr. \cite[Construction 3.10]{CDPI}}]\label{Constr1}
Using CW-complexes, I shall extend the previous construction to \emph{non-square-free monic monomials}. One proceeds as follows. Let $g=x_1^{j_1}\cdot\dotsc\cdot x_m^{j_m}$ be a degree $d=j_1+\dotsc+j_m$ monomial. Consider the following finite set ${V=\left\{x^1_1,\dotsc,x_m^{j_m}\right\}}$, let $\Delta(g)=\sP(V)$ be the associated abstract (finite) simplex, one considers the corresponding (topological) simplex (which is a CW-complex) $\widetilde{\Delta(g)}$.\smallskip

\noindent For any $k\in\{1,\dotsc,m\}$, if $j_k\leq1$ one does nothing, while if $j_k\geq2$, one recursively identifies, for $p$ varying from $0$ to $j_k-2$, the $p$-faces of the subsimplex ${\widetilde{\sP\left(\left\{x^1_k,\dotsc,x^{j_k}_k\right\}\right)}\subseteq\widetilde{\Delta(g)}}$ and the $q$-faces of $\widetilde{\Delta(g)}$ whose subfaces are identified: start with $p=0$, and one identifies all the $j_k$ points to one point, call it $x_k$. Then, for $p=1$, one obtains a bouquet of $\displaystyle\binom{j_k+1}{2}$ circles, and one identifies them in just one circle $\bS^1$ passing through $x_k$; and for $q=1$ one identifies all path with endpoints $x_h$ and $x_k$, where $h\neq k\in\{1,\dotsc,n\}$; and so on, up to the facets of $\widetilde{\sP\left(\left\{x^1_k,\dotsc,x^{j_k}_k\right\}\right)}$, \emph{i.e.} its $j_k+1$ $(j_k-1)$-faces, which, by construction, have all their boundaries in common, and one identifies all of them. Of course, one continues to identify the $q$-faces of $\widetilde{\Delta(g)}$ as indicated above.\smallskip

\noindent In this way, one obtains a finite CW-complex $\zeta_g$ of dimension $d-1$ obtained from the $(d-1)$-dimensional simplex $\widetilde{\Delta(g)}$, which $0$-skeleton $\zeta_g^0$ is $\left\{x_i\mid j_i\neq0\right\}$.\smallskip

\noindent I clarify this construction via some example.
\begin{example}
Let $g=x^2y$, then $\Delta(g)=\sP\left(\left\{x,x^2,y\right\}\right)$. I start identifying $\{x\}$ and $\left\{x^2\right\}$. And I finish to identify $\left\{x,y\right\}$ and $\left\{x^2,y\right\}$.
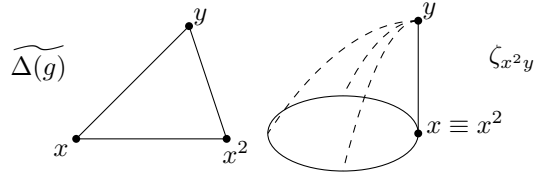
\begin{center}
\begin{tikzpicture}[scale=0.50]
\path (0,0) edge[left] (4,0);
\path (4,0) edge[right] (3,3);
\path (0,0) edge[left] (3,3);
\fill (0,0) circle(3pt);
\fill (4,0) circle(3pt);
\fill (3,3) circle(3pt);
\node at (-1,2){$\widetilde{\Delta\left(x^2y\right)}$};
\node at (-0.4,-0.3){$x$};
\node at (4.25,-0.35){$x^2$};
\node at (3.3,3.25) {$y$};
\end{tikzpicture}
\begin{tikzpicture}[scale=0.50]
\draw (-1,0) ellipse (2cm and 1cm);
\draw[dashed] (1,3) parabola (-3,0);
\draw[dashed] (1,3) parabola (-1,1);
\draw[dashed] (1,3) parabola (-1,-1);
\path (1,0) edge[left] (1,3);
\fill (1,0) circle(3pt);
\fill (1,3) circle(3pt);
\node at (3.5,2){$\zeta_{x^2y}$};
\node at (2.25,0.3){$x\equiv x^2$};
\node at (1.3,3.25){$y$};
\end{tikzpicture}
\end{center}
\end{example}
\begin{example}
Let $g=x^2y^2$, then $\Delta(g)=\sP\left(\left\{x,x^2,y,y^2\right\}\right)$.
\begin{center}
\begin{tikzpicture}[scale=0.50]
\path (0,0) edge[dashed,left] (6,2.5);
\path (0,0) edge[left] (4,0);
\path (4,0) edge[right] (6,2.5);
\path (4,0) edge[right] (3,4);
\path (3,4) edge[left] (6,2.5);
\path (0,0) edge[left] (3,4);
\fill (0,0) circle(3pt);
\fill (4,0) circle(3pt);
\fill (6,2.5) circle(3pt);
\fill (3,4) circle(3pt);
\node at (-1,2){$\widetilde{\Delta(g)}$};
\node at (-0.4,-0.3){$x$};
\node at (4.25,-0.35){$x^2$};
\node at (3.3,4.25) {$y$};
\node at (6.35,2.25){$y^2$};
\end{tikzpicture}
\end{center}
I start identifying $\{x\}$ and $\left\{x^2\right\}$, $\{y\}$ and $\left\{y^2\right\}$. At second step, I identify $\left\{x,y\right\},\left\{x^2,y\right\}$, $\left\{x,y^2\right\}$ and $\left\{x^2,y^2\right\}$. At the end I identify $\left\{x,x^2,y\right\}$ and $\left\{x,x^2,y^2\right\}$, $\left\{x,y,y^2\right\}$ and $\left\{x^2,y,y^2\right\}$.
\end{example}
\begin{example}
Let $g=x^3y$, then $\Delta(g)=\sP\left(\left\{x,x^2,x^3,y\right\}\right)$.
\begin{center}
\begin{tikzpicture}[scale=0.50]
\path (0,0) edge[dashed,left] (6,2.5);
\path (0,0) edge[left] (4,0);
\path (4,0) edge[right] (6,2.5);
\path (4,0) edge[right] (3,4);
\path (3,4) edge[left] (6,2.5);
\path (0,0) edge[left] (3,4);
\fill (0,0) circle(3pt);
\fill (4,0) circle(3pt);
\fill (6,2.5) circle(3pt);
\fill (3,4) circle(3pt);
\node at (-1,2){$\widetilde{\Delta(g)}$};
\node at (-0.4,-0.3){$x$};
\node at (4.25,-0.35){$x^2$};
\node at (3.3,4.5){$x^3$};
\node at (6.35,2.25){$y$};
\end{tikzpicture}
\end{center}
I start identifying $\{x\}$, $\left\{x^2\right\}$ and $\left\{x^3\right\}$. At second step, I identify $\left\{x,x^2\right\}$, $\left\{x,x^3\right\}$ and $\left\{x^2,x^3\right\}$. At third step I identify $\left\{x,y\right\}$, $\left\{x^2,y\right\}$ and $\left\{x^3,y\right\}$. At the end I identify $\left\{x,x^2,y\right\}$, $\left\{x,x^3,y\right\}$ and $\left\{x^2,x^3,y\right\}$.
\end{example}
\begin{example}
Let $g=x^2yz$, then $\Delta(g)=\sP\left(\left\{x,x^2,y,z\right\}\right)$.
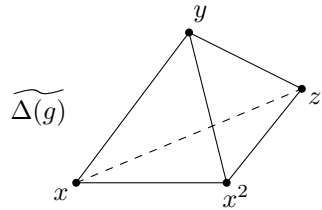
\begin{center}
\begin{tikzpicture}[scale=0.50]
\path (0,0) edge[dashed,left] (6,2.5);
\path (0,0) edge[left] (4,0);
\path (4,0) edge[right] (6,2.5);
\path (4,0) edge[right] (3,4);
\path (3,4) edge[left] (6,2.5);
\path (0,0) edge[left] (3,4);
\fill (0,0) circle(3pt);
\fill (4,0) circle(3pt);
\fill (6,2.5) circle(3pt);
\fill (3,4) circle(3pt);
\node at (-1,2){$\widetilde{\Delta(g)}$};
\node at (-0.4,-0.3){$x$};
\node at (4.25,-0.35){$x^2$};
\node at (3.3,4.5){$y$};
\node at (6.35,2.25){$z$};
\end{tikzpicture}
\end{center}
I start identifying $\{x\}$ and $\left\{x^2\right\}$. At second step, I identify $\left\{x,y\right\}$ and $\left\{x^2,y\right\}$, $\left\{x,z\right\}$ and $\left\{x^2,z\right\}$. At the end I identify $\left\{x,y,z\right\}$ and $\left\{x^2,y,z\right\}$.
\end{example}
\begin{remark}
For any index $k$, under this identification, the closure of each $(j_k-1)$-cell $\overline{\left\{x_k^1,\dotsc,x_k^{j_k}\right\}}$ becomes a point if $j_k=1$, a circle $\bS^1$ if $j_k=2$, a topological space with fundamental group $\Z_3$ (in particular, this is not a topological surface) if $j_k=3$, \emph{e.o.}
\end{remark}
\begin{proposition}[{cfr. \cite[Proposition 3.11]{CDPI}}]\label{prop4.1}
Every power in $x_k^{j_k}$ corresponds to a $\zeta_{x_k^{j_k}}$, and \emph{vice versa}.
\end{proposition}
\begin{corollary}
Let $\zeta_{x^d}$ be the CW-complex associated to monomial $x^d$. Then for any $d\in\N_{\geq1}$, $k\in\{0,\dotsc,d-1\}$, $\#\zeta_{x^d}^k=k+1$.
\end{corollary}
\noindent More in general:
\begin{proposition}[{cfr. \cite[Proposition 3.14]{CDPI}}]\label{prop4.2}
$\zeta_{x_1^{h_1}\cdot\dotsc\cdot x_m^{h_m}}\subseteq\zeta_{x_1^{k_1}\cdot\dotsc\cdot x_m^{k_m}}$ if and only if $x_1^{h_1}\cdot\dotsc\cdot x_m^{h_m}$ divides $x_1^{k_1}\cdot\dotsc\cdot x_m^{k_m}$.
\end{proposition}
\begin{prf}
Let $n=1$, this is Proposition \ref{prop4.1}. Let assume $n\geq2$. If ${x_1^{h_1}\cdot\dotsc\cdot x_m^{h_m}}$ divides $x_1^{k_1}\cdot\dotsc\cdot x_m^{k_m}$ then $\sP\left(\left\{x_1,\dotsc,x_m^{h_n}\right\}\right)\subseteq\sP\left(\left\{x_1,\dotsc,x_m^{k_n}\right\}\right)$; by the previous construction $\zeta_{x_1^{h_1}\cdot\dotsc\cdot x_m^{h_m}}\subseteq\zeta_{x_1^{k_1}\cdot\dotsc\cdot x_m^{k_m}}$. \emph{Vice versa}, if $\zeta_{x_1^{h_1}\cdot\dotsc\cdot x_m^{h_m}}\subseteq\zeta_{x_1^{k_1}\cdot\dotsc\cdot x_m^{k_m}}$ then ${\zeta_{x_1^{h_1}},\dotsc,}$ ${\zeta_{x_m^{h_m}}\subseteq\zeta_{x_1^{h_1}\cdot\dotsc\cdot x_m^{h_m}}\subseteq\zeta_{x_1^{k_1}\cdot\dotsc\cdot x_m^{k_m}}}$ hence $\zeta_{x_{\ell}^{h_{\ell}}}\subseteq\zeta_{x_{\ell}^{k_{\ell}}}$ for any $\ell\in\{1,\dotsc,m\}$. By the previous reasoning one has the claim.
\end{prf}
\noindent It is clear how to glue two of these finite CW-complexes: let $X=\zeta_{x_1^{j_1}\cdot\dotsc\cdot x_m^{j_m}}$ and $Y=\zeta_{x_1^{k_1}\cdot\dotsc\cdot x_m^{k_m}}$, of degree $j_1+\dotsc+j_m$ and $k_1+\dotsc+k_m$, one glues $X$ and $Y$ along the CW-complex $Z$ associated to lowest common divisor of $x_1^{j_1}\cdot\dotsc\cdot x_m^{j_m}$ and $x_1^{k_1}\cdot\dotsc\cdot x_m^{k_m}$.
\begin{remark}
If $x_1^{j_1}\cdot\dotsc\cdot x_m^{j_m}$ and $x_1^{k_1}\cdot\dotsc\cdot x_m^{k_m}$ are coprime monomials then one assumes $Z=\emptyset$.
\end{remark}
\noindent Finally, taking all these finite CW-complexes together, one obtains a CW-complex $P(m)$ in the following way:
\begin{align*}
C&=\bigsqcup_{x_1^{j_1}\cdot\dotsc\cdot x_m^{j_m}\in R}\zeta_{x_1^{j_1}\cdot\dotsc\cdot x_m^{j_m}} & P(m)={C}/{\sim}
\end{align*}
where $\sim$ is the equivalence relation induced by the above gluing.
\begin{proposition}[{cfr. \cite[Proposition 3.13]{CDPI}}]
There is bijection between the monomials of degree $d$ in $R$ and the $(d-1)$-cells of $P(m)$.
\end{proposition}
\noindent In other words, if one defines
\begin{displaymath}
\rho^{\prime}_d=\left\{f\in R_d\mid f\neq0\text{\, is a monic monomial}\right\}
\end{displaymath}
one has a bijection
\begin{align*}
\sigma^{\prime}_d\colon\rho^{\prime}_d &\to P(m)^{d-1}\setminus P(m)^{d-2}\\
x_1^{j_1}\cdot\dotsc\cdot x_m^{j_m}&\mapsto \zeta_{x_1^{j_1}\cdot\dotsc\cdot x_m^{j_m}}.
\end{align*}
\end{construction}
\noindent Let $\displaystyle f=\sum_{i=0}^Na_ig_i\in R_d$ be a degree $d$ homogeneous polynomial; for clarity, $g_i$ is a monomial of degree $d$ and $a_i\neq0$ for any $i$. One can associate to $f$ a finite $(d-1)$-dimensional CW-complex $\zeta_f$, where the $(d-1)$-skeleton is given by the $\zeta_{g_i}$'s glued together with the above procedure.
\begin{remarks}
\,\begin{enumerate}[a)]
\item The elements of $(d-1)$-skeleton of $\zeta_f$ correspond to elements of $(d-1)$-skeleton of $P(m)$.
\item Construction \ref{Constr1} generalizes the analogous one given in \cite[Section 3]{SR:1} and \cite[Section 1]{R:GA}. \hfill{$\Diamond$}
\end{enumerate}
\end{remarks}

\subsection{CW-complexes and annihilator ideals}

\begin{remark}\label{rem4.1}
In order to state the next theorem, I observe that the canonical basis of $R_d$ and $Q_d$ given by monic monomials are dual basis each other, \emph{i.e.}
\begin{displaymath}
X_1^{j_1}\cdot\dotsc\cdot X_m^{j_m}\left(x_1^{i_1}\cdot\dotsc\cdot x_m^{j_m}\right)=\delta^{i_1,\dotsc,i_m}_{j_1,\dotsc,j_m}
\end{displaymath}
where $i_1+\cdot\dotsc\cdot+i_m=j_1+\cdot\dotsc\cdot+j_m=d$ and $\delta^{i_1,\dotsc,i_m}_{j_1,\dotsc,j_m}$ is the Kronecker delta.\smallskip

\noindent This simple observation allows one to identify the \emph{dual differential operator $G_r$} of the monic monomial $g_r$ with the same element of the $(d-1)$-skeleton of $\zeta_f$ associated to $g_r$. To be clear, given a homogeneous polynomial $\displaystyle{f=\sum_{r=1}^ng_r\in R_d}$, monomial $G_r\in Q_d$ is such that $G_r(\lambda g_r)=\lambda$ for some $\lambda\in\K\setminus\{0\}$, and $G_r(g)=0$ for any other monomial $g\in R_d\setminus\left\langle g_r\right\rangle$. In other words, one associates to $g_r=x_1^{i_1}\cdot\dotsc\cdot x_m^{j_m}$ and to $G_r=X_1^{i_1}\cdot\dotsc\cdot X_m^{i_m}$ the CW-subcomplex $\zeta_{x_1^{i_1}\cdot\dotsc\cdot x_m^{i_m}}$ of $\zeta_f$.
\end{remark}
\begin{theorem}\label{th4.1}
Let $M_1,\dotsc,M_{\tau(m,d)}$ be the monic monomial basis of $R_d$, let $\displaystyle{f=\sum_{i\in I}M_i\in R_d}$ be a degree $d$ homogenous polynomial, where $I\subseteq\{1,\dotsc,\tau(n,d)\}$. Let $\zeta_f$ be the CW-complex associated to $f$. Then $I=\Ann(f)$ is generated by:
\begin{enumerate}[a)]
\item\label{th4.1.a} $X_k^{d+1}$, for any $k\in\{1,\dotsc,m\}$;
\item\label{th4.1.b} the monomials $X_1^{s_1}\cdot\dotsc\cdot X_m^{s_m}$ such that $s_1+\dotsc+s_m=j$, where $\zeta_{x_1^{s_1}\cdot\dotsc\cdot x_m^{s_m}}$ is a (minimal) element of the $(j-1)$-skeleton of $P(m)$ not contained in $\zeta_f$ (for any $j\in\{1,\dotsc,d\}$);
\item\label{th4.1.c} the degree $j$ homogeneous polynomials $P_1-P_2$ with all coefficients equal to $1$ such that $\zeta_{P_1(f)}=\zeta_{P_2(f)}\subseteq\zeta_f$.
\end{enumerate}
\end{theorem}
\begin{remarks}
About the polynomials described in the previous item \ref{th4.1.c}, I highlight the following facts:
\begin{enumerate}[a)]
\item they can not be coprime. For example, let $m=d=2$ and $f=x_1^2+x_1x_2$ then $\Ann(f)=\left(X_1^3,X_2^2,X_1^2-X_1X_2\right)$.
\item they can not be binomials. For example, let $m=4$, $d=3$ and $\displaystyle f=\sum_{1\leq i<j<k\leq4}x_ix_jx_k$ then
\begin{displaymath}
\Ann(f)\cap\K\left[X_1,X_2,X_3,X_4\right]_2=\left\langle X_iX_j+X_kX_{\ell}-X_iX_{\ell}-X_jX_k\mid\{i,j,k,\ell\}=\{1,2,3,4\}\right\rangle
\end{displaymath}
but for each $\{i,j,k,\ell\}=\{1,2,3,4\},X_iX_j-X_jX_k\notin\Ann(f)$.
\end{enumerate}
\noindent Both the previous examples determine SGAG algebras. \hfill $\Diamond$
\end{remarks}

\noindent \textbf{Proof of Theorem \ref{th4.1}}. Let $A=Q/I$. Since $A$ has socle degree $d$, then $X_k^{d+1}\in I$ for any ${k\in\{0,\dotsc,m\}}$.\smallskip

\noindent Using the identification introduced in Remark \ref{rem4.1}, the monomials $X_1^{s_1}\dotsc X_m^{s_m}$ such that $s_1+\dotsm+s_m=j$ and $\zeta_{x_1^{s_1}\dotsc x_m^{s_m}}$ is a (minimal) element of the $(j-1)$-skeleton of $P(m)$ not contained in $\zeta_f$ belong to $I$.\smallskip

\noindent Let $\alpha$ be the homogeneous ideal generated by the monomials given by items \ref{th4.1.a} and \ref{th4.1.b}. One has $\displaystyle A\cong\frac{Q}{\alpha}/\frac{I}{\alpha}$ and the short exact sequence
\begin{displaymath}
\xymatrix{
0\ar[r] & \displaystyle\frac{I}{\alpha}\ar[r] & \displaystyle\frac{Q}{\alpha}\ar[r] & A\ar[r] & 0
}
\end{displaymath}
given by evaluation at $f$, \emph{i.e.} the residual class of $Y\in Q$ modulo $\alpha$ is sent to ${Y(f)\in A}$. Let $P\in Q$ be a polynomial with all coefficients equal to $1$, then it is different from $0$ modulo $\alpha$ if and only if $P(f)\in\K\setminus\{0\}$ or $\zeta_{P(f)}\cap\zeta_f\neq\emptyset$. Without loss of generality, one assumes $\zeta_{P(f)}\subseteq\zeta_f$. Thus let $P_1,P_2\in Q_j$ be polynomials with all coefficients equal to $1$, where $j\in\{1,\dotsc,d\}$, such that $\zeta_{P_1(f)},\zeta_{P_2(f)}\subseteq\zeta_f$ hence $P_1-P_2\in I$ if and only if $P_1(f)=P_2(f)$. By hypothesis this condition is equivalent to $\zeta_{P_1(f)}=\zeta_{P_2(f)}\subseteq\zeta_f$ \emph{i.e.} the claim holds. \hfill $\Box$

\section{On the Hilbert vector of GAG algebras}
\markboth{On th Hilbert vector of GAG algebras}{On th Hilbert vector of GAG algebras}

\noindent Consider the family $\cAG_{\K}(m,d)$ of GAG $\K$-algebras of socle degree $d$ and codimension $m$. Consider the family of symmetric length $d+1$ vectors of type ${\left(1,h_1,\dotsc,h_1,1\right)\in\left(\N_{\geq1}\right)^{d+1}}$. There is a natural partial order on the set of these vectors
\begin{displaymath}
\left(1,h_1,\dotsc,h_1,1\right)\preceq\left(1,\wh_1,\dotsc,\wh_1,1\right)
\end{displaymath}
if $h_i\leq\wh_i$ for any $\displaystyle i\in\left\{2,\dotsc,\left\lfloor\frac{d}{2}\right\rfloor\right\}$. This order can be restricted to $\cAG_{st,\K}(m,d)$ which becomes a \emph{poset}.
\begin{definition}
Let $m,d\in\N_{\geq0}$ be fixed and let $H$ be a length $d+1$ symmetric vector. One says $H$ is \emph{the minimum} (\emph{a minimal}, respectively) Artinian Gorenstein Hilbert vector of socle degree $d$ and codimension $m$ if there is a GAG algebra such that $\Hilb(A)=H$ and $H$ is the minimum (minimal, respectively) in $\cAG_{\K}(m,d)$ with respect to $\preceq$. (To be precise, if $\wH$ is a comparable Artinian Gorenstein Hilbert vector such that $\wH\preceq H$, then $\wH=H$.)
\end{definition}
\begin{lemma}\label{lem5.4}
Let $\ua=(a_i)\in\left(\K\setminus\{0\}\right)^I$, let $M_1,\dotsc,M_{\tau(m,d)}$ be the monic monomial basis of $R_d$, let $\displaystyle f_{\ua}=\sum_{i\in I}a_iM_i\in R_d$ be a degree $d$ homogenous polynomial, where $I\subseteq\{1,\dotsc,\tau(m,d)\}$, let $f=f_{(1,\dotsc,1)}$ for short. Then the Hilbert vector of $Q/\Ann(f)$ is the minimum between the Hilbert vectors of $Q/\Ann\left(f_{\ua}\right)$'s.
\end{lemma}
\begin{prf}
For sake of simplicity, let $\Ann(f)=I$ and $\Ann\left(f_{\ua}\right)=I_{\ua}$. Let $\zeta_f$ be the CW-complex associated to $f$.\smallskip

\noindent By Theorem \ref{th4.1}, the monomials described by items (\ref{th4.1.a}) and (\ref{th4.1.b}) belongs to $I_{\ua}$.\smallskip

\noindent Let $\displaystyle\wP_1=\sum_{k=1}^{s_1}b_{1,k}\wG_{1,k},\wP_2=\sum_{\ell=1}^{s_2}b_{2,\ell}\wG_{2,\ell}\in Q_j$ be polynomials, where ${j\in\{1,\dotsc,d\}}$, such that $\wG_{1,k}$'s and $\wG_{2,\ell}$'s determine a subsystem of the monic monomial basis of $Q_j$, and $\wP_1\left(f_{\ua}\right)=\wP_2\left(f_{\ua}\right)$ \emph{i.e.} $\wP_1-\wP_2\in I_{\ua}$. Let $\displaystyle P_1=\sum_{k=1}^{s_1}\wG_{1,k}$, $\displaystyle{P_2=\sum_{\ell=1}^{s_2}\wG_{2,\ell}\in Q_j}$; without loss of generality, let $\zeta_{P_1(f)},\zeta_{P_2(f)}\subseteq\zeta_f$, \emph{i.e.} both $\wG_{1,k}$ and $\wG_{2,\ell}$ do not vanish neither on $f_{\ua}$ nor on $f$. Since $\wP_1\left(f_{\ua}\right)=\wP_2\left(f_{\ua}\right)$ then each monomial on left hand side is equal to a monomial on right hand side; forgetting the coefficient of any monomial one has $P_1(f)=P_2(f)$ that is $P_1-P_2\in\Ann(f)$.\smallskip

\noindent From all this, one has an injective function from $I_{\ua}$ to $I$ which sends each generator described by Theorem \ref{th4.1} to itself forgetting the relevant coefficients. Let ${H=\left(1,h_1,\dotsc,h_1,1\right)}=\Hilb(Q/I)$ and let $H_{\ua}=\left(1,h_{\ua,1},\dotsc,h_{\ua,1},1\right)=\Hilb\left(Q/I_{\ua}\right)$; then by the previous reasoning $h_i\leq h_{\ua,i}$ for any $\displaystyle i\in\left\{1,\dotsc,\left\lfloor\frac{d}{2}\right\rfloor\right\}$ \emph{i.e.} $H\preceq H_{\ua}$.
\end{prf}
\begin{remarks}
Without changing the previous notations:
\begin{enumerate}[a)]
\item one has for any $\lambda\in\K\setminus\{0\}$, $\Ann(f)=\Ann(\lambda f)$.
\item the function from $I_{\ua}$ to $I$ is not surjective, in general. For example, let $m=2$, ${d=3}$ and $f_{(1,3,3,1)}=x_1^3+3x_1^2x_2+3x_1x_2^2+x_2^3$ hence ${f=x_1^3+x_1^2x_2+x_1x_2^2+x_2^3}$. One has $x_1^2-x_2^2\in\Ann(f)$ but for each $\lambda_1,\lambda_2\in\K\setminus\{0\}$, ${\lambda_1x_1^2-\lambda_2x_2^2\notin\Ann\left(f_{(1,3,3,1)}\right)}$.~\hfill{$\Diamond$}
\end{enumerate}
\end{remarks}
\begin{proposition}
Let $M_1,\dotsc,M_{\tau(m,d)}$ be the monic monomial basis of $R_d$, let $\displaystyle f=\sum_{i\in I}M_i\in R_d$ be a degree $d$ homogenous polynomial, where $I\subseteq\{1,\dotsc,\tau(m,d)\}$. Let $\zeta_f$ be the CW-complex associated to $f$ and let $A=Q/\Ann(f)$. Then
\begin{displaymath}
A=\bigoplus_{h=0}^dA_h,\,1\leq\dim A_h\leq s_h
\end{displaymath}
where $s_h$ is the number of the $(h-1)$-cells of $\zeta_f$ (with the convention $s_0=1$).
\end{proposition}
\begin{prf}
Trivially $A_0\cong\K$. By definition, if $h\in\{1,\dotsc,d\}$, $A_h$ is generated by the (canonical images of the) monomials $\Omega_s\in Q_h$ that do not annihilate $f$. This means that, if one writes
\begin{displaymath}
\Omega_l=X_1^{\ell_1}\dotsm X_m^{\ell_m}\,\text{\,where\,}\,\ell_1+\dotsm+\ell_m=h,
\end{displaymath}
there exists $s\in\{1,\dotsc,m\}$ such that $g_s=x_1^{\ell_1}\dotsm x_m^{\ell_m}g^{\prime}_s$, where $g^{\prime}_s\in R_{d-h}$ is a (nonzero) monomial; this means that $X_1^{\ell_1}\dotsm X_m^{\ell_m}$ is an element of the $(h-1)$-skeleton of $\zeta_f$ by Proposition \ref{prop4.2}.
\end{prf}

\subsection{On the entries of the Hilbert vector of SGAG algebras}

\noindent Let $d\in\N_{\geq4}$, $n\in\N_{\geq2}$, let $\cAG_{st,\K}(n,d)$ be the set of SGAG $\K$-algebras of socle degree $d$ and codimension $m$. One knows that the Hilbert vector of $A\in\cAG_{st,\K}(n,d)$ is symmetric by Poincar\'e duality and the $0$-th and $1$-st entries are $1$ and $n$, respectively. Therefore it is matter to the aim of this paper to study the other entries of this vector.\smallskip

\noindent By Theorem \ref{th1.1} there exists $[f]\in U^{(n,d)}_{st}$ (see Equation \eqref{eq1}) such that ${A\cong Q/\Ann(f)}$. Consider the short exact sequence
\begin{displaymath}
\xymatrix{
0\ar[r] & \Ann(f)_k\ar[r] & Q_k\ar[r] & A_k\ar[r] & 0,
}
\end{displaymath}
where $k\in\{2,\dotsc,d-2\}$. It is
\begin{displaymath}
h_k=\dim A_k=\dim Q_k-\dim\Ann(f)_k,
\end{displaymath}
where I remind $\displaystyle\dim Q_d=\tau(n,d)=\binom{n+d-1}{d}$. From all this, in order to minimize $h_k$ one needs to maximize $\dim\Ann(f)_k$.
\begin{lemma}\label{lem5.1}
Let $k\in\{2,\dotsc,d-2\}$. Let $U$ be an irreducible component of $U^{(n,d)}_{st}$ (see Equation \eqref{eq1}), let $M_1,\dotsc,M_{\tau(n,d)}$ be the monic monomial basis of $\K\left[x_1,\dotsc,x_n\right]_d=R_d$ and let $[f]\in U$. Assume that
\begin{displaymath} 
f=\sum_{i\in I}a_iM_i\,\text{\,where\,}\,I\subseteq\{1,\dotsc,\tau(n,d)\}
\end{displaymath}
then
\begin{displaymath}
\dim\left\langle q(f)\in R_{d-k}\mid q\in Q_2\right\rangle\leq\dim\left\langle q\left(M_i\right)\in R_{d-k}\mid q\in Q_k,i\in I\right\rangle.
\end{displaymath}
Moreover, under opportune hypothesis, if $U=U^{(n,d)}_{FP}$ (the full Perazzo locus, see Equation \eqref{eq5}) then between the annihilators of the elements of $\K\left[x_1,\dotsc,x_{n+\tau(n,d-1)}\right]_d$, $\Ann(f)_k$ is generated by the maximum number of degree $k$ monomials.
\end{lemma}
\begin{prf}
The inequality between the dimensions is trivial by the bilinearity of
\begin{displaymath}
B_k\colon(q,p)\in Q_k\times R_d\to q(p)\in R_{d-k};
\end{displaymath}
in particular, the vector space on the right hand side of that inequality has minimal dimension between the vector subspaces of $R_{d-k}$ such that: they are generated by the images of degree $d$ monomials via $B_k$, and these monomials determine a polynomial $g$ such that $[g]\in U^{(m,d)}_{st}$.\smallskip

\noindent Assume that $U$ is the full Perazzo locus, if there exists a polynomial $g$ such that $[g]\in U^{(n+\tau(n,d-1),d)}_{st}$ and $\Ann(g)_k$ is generated by more degree $k$ monomials than $\Ann(f)_k$, then this means that $g$ is determined by less $\#I$ monomials than $f$. Indeed, by Lemma \ref{lem5.4}, let
\begin{displaymath}
g=\sum_{j\in J}M_j
\end{displaymath}
where $J\subseteq\{1,\dotsc,\tau(m,d)\}$ and $M_j$'s determine a subsystem of the monomial basis of $\K\left[x_1,\dotsc,x_{n+\tau(n,d-1)}\right]_d$ and let
\begin{displaymath}
L=\left\{\left(\ell_1,\dotsc,\ell_k\right)\in\{1,\dotsc,n+\tau(n,d-1)\}^k\mid\ell_1\leq\dotsc\leq\ell_k,x_{\ell_1}\cdot\dotsc\cdot x_{\ell_k}\in\Ann(g)_k\right\};
\end{displaymath}
then
\begin{displaymath}
\forall j\in J,\left(\ell_1,\dotsc,\ell_k\right)\in L,\,M_j=x^{\uj},\sum_{i=1}^k\pr_{\ell_i}\left(\uj\right)\leq k-1\,\text{\,where\,}\,\uj\in T\left(n+\tau(n,d-1),d\right);
\end{displaymath}
for the notations see Equations \eqref{eq2}, \eqref{eq4} and \eqref{eq3}. In this way, one has $\#J\lneqq\#I$ because the monomials indexed by $J$ have to satisfy more conditions given by $\pr$'s functions than monomials indexed by $I$. By Lemma \ref{lem2.2}, all this is impossible hence one has the claim.
\end{prf}

\subsection{Full Perazzo algebras case}

\noindent Moreover, the previous lemma can be improved as it follows.
\begin{lemma}\label{lem5.2}
Let $k\in\{2,\dotsc,d-2\}$. Let $f$ be a full Perazzo polynomial. Between the annihilators of the elements of $\K\left[x_1,\dotsc,x_{\tau(n,d-1)},u_1,\dotsc,u_n\right]_d$, $\Ann(f)_k$ is generated by the maximum number of degree $k$ monomials and binomials.
\end{lemma}
\begin{prf}
It is known
\begin{gather*}
\Ann(f)_k\supseteq\left\langle X_{i_1}U^{\uj_1}-X_{i_2}U^{\uj_2}\mid i_1,i_2\in\{1,\dotsc,\tau(n,d-1)\},\uj_1,\uj_2\in T(n,k-1),\right.\\
\left.x_{i_1}u^{\uj_1}P_{i_1,i_2}=x_{i_2}u^{\uj_2}P_{i_1,i_2},\,\text{\,where\,} P_{i_1,i_2}=\gcd\left(x_{i_1}M_{i_1},x_{i_2}M_{i_2}\right),\,\deg P_{i_1,i_2}=d-k\right\rangle\\
\text{\,and\,}\,M_1,\dotsc,M_{\tau(n,d-1)}\,\text{\,is the monic monomial basis of\,}\,\K\left[u_1,\dotsc,u_n\right]_{d-1},
\end{gather*}
and by \cite[Theorem 3.18]{CDPI}, these are all degree $k$ binomial generators of $\Ann(f)$.\smallskip

\noindent Let $U$ be another irreducible component of $U^{(n+\tau(n,d-1),d)}_{st}$ which dimension is $\tau(n,d-1)-1$, let $[g]\in U$; by Lemma \ref{lem5.1} $\Ann(g)_k$ is generated by the same number of degree $k$ monomials of $\Ann(f)_k$. Explicitly, without loss of generality (see Lemma \ref{lem5.4}) let
\begin{displaymath}
g=\sum_{\ell\in L}N_{\ell}
\end{displaymath}
where $L\subseteq\{1,\dotsc,\tau(m,d)\}$ and $N_{\ell}$'s form the monomial basis of $\K\left[x_1,\dotsc,x_{\tau(n,d-1)},u_1,\dotsc,u_n\right]_d$. Let $I$ the set of tuples $\left(i_1,i_2,\uj_1,\uj_2\right)$ such that $X_{i_1}U^{\uj_1}-X_{i_2}U^{\uj_2}\in\Ann(f)_k$, and let $J$ be the set of tuples $\left(\ua,\ub,\uc,\ud\right)$ such that $X^{\ua}U^{\uc}-X^{\ub}U^{\ud}\in\Ann(g)_k$. If $\#I\lneqq\#J$ then by the conditions
\begin{gather*}
\forall\ell\in L,\left(\ua,\ub,\uc,\ud\right)\in J,\,N_{\ell}=(xu)^{\ul},\sum_{\ul\in\wL}(xu)^{D_{\ua}D_{\uc}\left(\ul\right)}=\sum_{\ul\in\wL}(xu)^{D_{\ub}D_{\ud}\left(\ul\right)}\\
\text{\,with\,}\,\wL\subseteq T\left(n+\tau(n,d-1),d\right),
\end{gather*}
where for the notations see Equations \eqref{eq2}, \eqref{eq3} and $D_{\ua}D_{\uc}\left(\ul\right)$ ($D_{\ua}D_{\uc}\left(\ul\right)$, respectively) means that the indexes of $\ul$ which correspond to positions indicated in the vector $\left(\ua,\uc\right)$ ($\left(\ub,\ud\right)$, respectively) decrease of one unit (cfr. Equation \eqref{eq6}). From all this:
\begin{gather}\label{eq7}
\forall\left(\ua,\ub,\uc,\ud\right)\in J,\,\wL\,\text{\,is partitioned in pairs\,}\,\left\{\ul_{\ua,\uc},\ul_{\ub,\ud}\right\}\mid\\
\text{for each of such pair\,}\,D_{\ua}D_{\uc}\left(\ul_{\ua,\uc}\right)=D_{\ub}D_{\ud}\left(\ul_{\ub,\ud}\right) \notag
\end{gather}
hence one deduces $\#L=\#\wL\lneqq\tau(n,d-1)$, because the monomials indexed by $L$ have to satisfy more conditions than monomials determining a full Perazzo polynomial. All this is impossible, by the assumption on the dimension of $U$, therefore the claim follows.
\end{prf}
\noindent Now I am in position to state the main lemma of this subsection.
\begin{lemma}\label{lem5.3}
Let $k\in\{2,\dotsc,d-2\}$. Let $f$ be a degree $d$ full Perazzo polynomial. Then between the annihilators of the elements of $\K\left[x_1,\dotsc,x_{\tau(n,d-1)},u_1,\dotsc,u_n\right]_d$, $\Ann(f)_k$ is generated by the maximum number of degree $k$ polynomials.
\end{lemma}
\begin{prf}
By \cite[Theorem 3.18]{CDPI}, the degree $k$ monomials and binomials described in Lemmata \ref{lem5.1} and \ref{lem5.2} are all and the only generators of $\Ann(f)_k$.\medskip

\noindent Let assume there exists a polynomial $g$ such that $[g]\in U^{(n+\tau(n,d-1),d)}_{st}$ and $\Ann(g)_k$ is generated by more degree $k$ polynomials than $\Ann(f)_k$. By Lemma \ref{lem5.4}, let
\begin{displaymath}
g=\sum_{\ell\in L}M_{\ell}
\end{displaymath}
where $L\subseteq\{1,\dotsc,\tau(m,d)\}$ and $M_j$'s determine a subsystem of the monomial basis of $\K\left[x_1,\dotsc,x_{n+\tau(n,d-1)}\right]_d$; by Theorem \ref{th4.1}.\ref{th4.1.c} there exist $P_1,P_2\in\K\left[X_1,\dotsc,X_{\tau(n,d-1)},U_1,\dotsc,U_n\right]_k$ such that $P_1-P_2\in\Ann(g)_k$, in addition to other possible monomials and binomials. If $\dim\Ann(g)_k>\dim\Ann(f)_k$, \emph{id est} $g$ is annihilated by more degree $k$ polynomials than $f$ then this means that $g$ is determined by less monomials than $f$. Indeed, imitating the proof of the previous lemma, $L$ corresponds to a subset $\wL$ of $T\left(n+\tau(n,d-1),d\right)$; indexing the generators of $\Ann(f)_k$ and $\Ann(g)_k$ by sets $I$ and $J$, respectively, one suppose $\#I\lneqq\#J$ then $\wL$ is partitioned in tuples opportunely. This partition agrees with the generators of $\Ann(g)_k$ as explained in the Equation \eqref{eq7}. Thus $\#L=\#\wL\lneqq\tau(n,d-1)$, \emph{id est} $[g]$ belongs to an irreducible component of $U^{(n+\tau(n,d-1),d)}_{st}$ of dimension less than $U^{(n,d)}_{FP}$, and this is impossible by Lemma \ref{lem2.2}.
\end{prf}
\noindent After all this, I recall that in \cite{BGIZ} the authors posited the following \emph{Full Perazzo Conjecture}.
\begin{conjecture}[{\cite[Conjecture 2.6]{BGIZ}}]\label{conj}
Let $H$ be the Hilbert vector of a full Perazzo algebra of socle degree $d$ and codimension $n+\tau(n,d-1)$, where $d\in\N_{\geq4}$ and $n\in\N_{\geq2}$. Then $H$ is minimal in $\cAG_{st,\K}\left(n+\tau(n,d-1),d\right)$.
\end{conjecture}
\noindent By the reasoning at the begin of this section and by Lemma \ref{lem5.3}, the main theorem of this paper is proved.
\begin{theorem}\label{th5.1}
The Full Perazzo Conjecture holds.
\end{theorem}
\noindent The \emph{standard} assumption in the previous theorem cannot be skipped as the following example proves.
\begin{example}\label{ex2.1}
Let $n=2$ and $d=4$.
\begin{enumerate}
\item Let $f_1=u_1^2u_2^2$. It results
\begin{gather*}
\Ann\left(f_1\right)_1=\left\langle X_1,X_2X_3,X_4\right\rangle,\,\Ann\left(f_1\right)_2=\left\langle X_1^2,X_2^2,X_3^2,X_4^2,\right.\\
\left.X_1X_2,X_1X_3,X_1X_4,X_2X_3,X_2X_4,X_3X_4,X_1U_1,X_1U_2,X_2U_1,X_2U_2,X_3U_1,X_3U_2,X_4U_1,X_4U_2\right\rangle;
\end{gather*}
\noindent thus $H_1=\Hilb\left(\K\left[X_1,X_2,X_3,X_4,U_1,U_2\right]/\Ann\left(f_1\right)\right)=(1,2,3,2,1)$.
\item Let $f_2=x_1u_1^3+x_2u_1^2u_2+x_3u_1u_2^2+x_4u_2^3$. It results
\begin{gather*}
\Ann\left(f_2\right)_1=\{0\},\,\Ann\left(f_2\right)_2=\left\langle X_1^2,X_2^2,X_3^2,X_4^2,X_1X_2,X_1X_3,X_1X_4,X_2X_3,X_2X_4,X_3X_4,\right.\\
\left.X_1U_2,X_4U_1,X_1U_1-X_2U_2,X_2U_1-X_3U_2,X_3U_1-X_4U_2\right\rangle;
\end{gather*}
\noindent thus $H_2=\Hilb\left(\K\left[X_1,X_2,X_3,X_4,U_1,U_2\right]/\Ann\left(f_2\right)\right)=(1,6,5,6,1)$.
\end{enumerate}
By definition $H_1\precneqq H_2$ hence the previous statement holds.
\end{example}
\noindent By all this, Theorem \ref{th1.2} does not work for full Perazzo algebras.
\bigskip

\noindent{\bf Statement about competing or financial interests.} The author has no competing or financial interests to declare that are relevant to the content of this article.

\end{document}